\documentclass[a4paper,11pt]{article}
\usepackage{graphicx}
\usepackage[utf8]{inputenc}
\usepackage[T1]{fontenc}
\usepackage{lmodern}
\usepackage{xcolor}

\begin{document} 
\title{An inequality involving primes and the product of primes}
\author{{\small Andrej Leško}}
\date{}
\maketitle
\noindent
{\small
\textbf{Abstract:} In this note we prove an inequality involving primes and the product of consecutive primes.
\\
\textbf{Mathematics Subject Clasification:} 11A 41}
\\
\section{Introduction}
In article [\textbf{1}] from 2009 S. Zhang proved that \( p_{n+1}^{n-\pi(n)}>2^{p_{n+1}} \) for \( n\geq{20} \), and as a consequence of this result and Panaitopol's inequality [\textbf{2}] $ \prod_{k\leq{n}}p_k>p_{n+1}^{n-\pi(n)} $ for \( n\geq{2} \), deduced a lower bound for a product of consecutive primes, namely, \( \prod_{k\leq{n}}p_k>2^{p_{n+1}} \) for \( n\geq{10} \). As usual \( p_k \) denotes the k-th prime number, \( \pi(x) \) the prime counting function and \( \prod_{k\leq{n}}p_k \) the product of consecutive primes up to \( p_n \). We will introduce and prove a generalization of S. Zhang's results:
\\\\
\textbf{Theorem 1:} Let \( 1<c<e \) , (\textit{e=2.718-Euler number}) be a real number, then for every integer \( k\geq{0} \), there is a number \( N_k \) such that:
\\
\begin{equation}
 n^{n-\pi(n)}>c^{p_{n+k}}\textrm{, for }n>N_k .
\label{eq:main1}
\end{equation}
\\
Considering that \( p_{n+1}^{n-\pi(n)}>n^{n-\pi(n)} \) and taking \( k=1 \), \( c=2 \) we get S. Zhang's inequality, provided we independently confirm the result for numbers smaller than \( N_k \). 
\\
If we combine Theorem 1 with the Panaitopol's inequality, we get a new lower bound for a product of consecutive primes:
\\
\\
\textbf{Corollary 1:} Let \( 1<c<e \), be a real number and \( k\geq{0} \) positive integer, then there is an \( N_k \) such that 
\\
\begin{equation}
\prod_{j\leq{n}}p_j>c^{p_{n+k}} \textrm{, for }n\geq{N_k}. 
\end{equation}
\section{Proof of the main results}
\textbf{Proof of Theorem 1}
\\\\
The function
\begin{equation}
f_k(n)=1-\frac{\log{c}}{\log{n}}\left( 1+\frac{k}{n}\right) \log\left((n+k)\log(n+k)\right)-\frac{1.25506}{\log{n}}, n\geq{2}
\label{eq:th0}
\end{equation}
$(1.25506\approx{30\cdotp\frac{\log{113}}{113}})\quad$
has a limit \(\lim_{n\rightarrow\infty} f_k(n)=1-\log{c}>0 \) independent of \( k \), so there is \( N_k\geq{2} \) such, that \( f_k(n)>0 \) for \( n>N_k \). So we have
\\
\begin{equation}
 \frac{1.25506}{\log{n}}<1-\frac{\log{2}}{\log{n}}\left( 1+\frac{k}{n}\right) \log((n+k)\log(n+k)).
 \label{eq:th1}
\end{equation}
It is known that [\textbf{3}] 
\begin{equation}
\frac{1.25506}{\log{n}}>\frac{\pi(n)}{n}.    
\end{equation} 
After inserting this inequality into Equation~(\ref{eq:th1}) we find
\begin{equation}
	n-\pi(n)>\frac{\log{c}}{\log{n}}(n+k)\log((n+k)\log(n+k)).
\end{equation}
The well known estimate for the $n$-th consecutive prime ( $p_n<n\log(n\log{n})$ [\textbf{3}]) gives
\begin{equation}
	(n+k)\log((n+k)\log(n+k))>p_{n+k} \label{eq:th12}.
\end{equation}
Together with Equation~(\ref{eq:th12}) we have
\begin{equation}
	n-\pi(n)>\frac{\log{c}}{\log{n}}p_{n+k}, \quad\textrm{and finally}\qquad n^{n-\pi(n)}>c^{p_{n+k}}. 
\end{equation}
This completes the proof of \textbf{Theorem 1}.
\\
\\
Remark: Function \( f_k(x) \) used in the proof of the Theorem 1 allows us to estimate the number \( N_k \). Since \( f_k(x) \) is  continuous on the interval \( \left[ 2,\infty\right)  \) and increasing with \( f_k(2)<0 \) and \( \lim_{x\rightarrow\infty}{ f_k(x)}=1-\log{c}>0\) for \( 1<c<e \), this implies that \( f_k(x) \) has exactly one zero on that interval. Denoting this zero as $x_k$, we have from (\ref{eq:th0}), \(f_k(N_k)>0\) for \(N_k>x_k\). So we can take \(N_k\) equal to the integer part of \(x_k\).
This number need not be the smaller one  such that (\ref{eq:main1}) is valid. To find the lowest bound for \( N_k \) we must directly check numbers smaller than \( N_k \).
\\ \\
\textbf{Proof of Corollary 1}
\\
By Panaitopol ineq. \( \mathbf{\left[ 2\right]}  \) we have \(\prod_{j\leq{n}}p_j>p_{n+1}^{n-\pi(n)} \) , also \( p_{n+1}^{n-\pi(n)} >n^{n-\pi(n)}\), and finally by Theorem 1 \( \prod_{j\leq{n}}p_j >n^{n-\pi(n)}>c^{p_{n+k}}\) for \( n>N_k \).
\\
\section{Appendix}
 Let us mention that using the function \( f(x)=1-\log{2}\cdotp\left( 1+\frac{1}{x}\right) -\frac{1.71678}{\log(x\log{x})} \), we can directly prove original Zhang's result in the form \( p_{n+1}^{n-\pi(n)}>2^{p_{n+1}} \), for \( n\geq{20} \). Analysis shows that \( f(x) \) is increasing on interval \( \left[ 2,\infty\right)  \) with only zero at \( x\approx{74.39} \), so \( \frac{1.71678}{\log(n\log{n})}<1-\log{2}\cdotp\left( 1+\frac{1}{n}\right)  \) for \( n>74 \). It can be shown that \( \pi(n)<1.71678.\frac{n}{\log(n\log{n})} \) for  \(n\geq{2}  \) \textbf{(a.)}, so \( \frac{\pi(n)}{n}<1-\log{2}\cdotp\left( 1+\frac{1}{n}\right)  \). After multiplying by \( n \) and inserting \( \frac{\log((n+1)\log(n+1))}{\log((n+1)\log(n+1))} \) on\\\\ the right hand side, we can write \( n-\pi(n)>(n+1)\log{2}\frac{\log((n+1)\log(n+1))}{\log((n+1)\log(n+1))}\). From the well known estimates \( (n+1)\log((n+1)\log(n+1))>p_{n+1} \)\qquad and \( \log((n+1)\log(n+1))<p_{n+1} \) , it follows that \( (n-\pi(n))\log{p_{n+1}}>p_{n+1}\log{2} \) or \( p_{n+1}^{n-\pi(n)}>2^{p_{n+1}} \) . We have only to manualy confirm the identity for numbers smaller than 71.\\\\
\textbf{Proof of inequality (a.)}: For \( n\geq{2} \) we have \( \pi(n)<1.25506\frac{n}{\log{n}} \mathbf{\left[  3\right] } \), thus, \\
\( \pi(n)\frac{\log(n\log{n})}{n}<1.25506\frac{n}{\log{n}}\cdotp\frac{\log(n\log{n})}{n} =1.25506\cdotp\left( 1+\frac{\log\log{n}}{\log{n}}\right) \). Function \( 1+\frac{\log\log{x}}{\log{x}} \) has max. value \( 1+\frac{1}{e} \) at \( e^e \), so\\ \( \pi(n)\frac{\log(n\log{n})}{n}<1.25506\cdotp\left( 1+\frac{1}{e}\right) \approx{1.71678} \).
\\\\\\\\\\\\\\
\textbf{References}
\\
\( \mathbf{\left[ 1\right] } \) Shaohua Zhang, Anew inequality involving primes, arXiv:0908.2943v1, (2009)	
\\
\( \mathbf{\left[ 2\right] } \) Panaitopol Laurentiu, An inequality involving prime numbers, Univ. Beograd, Publ. Elektrotehn. Fak. Ser. Mat.,11,33-35 (2000)
\\
\( \mathbf{\left[ 3\right] } \) Rosser J. Barkley, Shoenfeld Lowell, Approximate formulas for some functions of prime numbers, Illinois J.Math.,6,64-94, (1962)
\end{document}